\documentclass[reqno]{amsart}    
\usepackage{enumerate,amsmath}    
\usepackage{pb-diagram} 
\usepackage{pstricks,pst-node,pst-coil,pst-plot}  
\theoremstyle{plain}     
\newtheorem{step}{Step}
\newtheorem{thm}{Theorem}[section]    
\newtheorem{theorem}[thm]{Theorem}    
\newtheorem{cor}[thm]{Corollary}    
\newtheorem{corollary}[thm]{Corollary}    
\newtheorem{lem}[thm]{Lemma}    
\newtheorem{lemma}[thm]{Lemma}    
    
\newtheorem{proposition}[thm]{Proposition}    
    
\theoremstyle{remark}     
    
\newtheorem{remark}[thm]{Remark}    
\newtheorem{examples}[thm]{Examples}    
\theoremstyle{definition}     
    
\newtheorem{definition}[thm]{Definition}    
\def\proof#1{{\par\medbreak\noindent {\it Proof\setbox0\hbox{#1}%
\ifdim\wd0=0pt .\else\ \ignorespaces #1.\fi}\enspace}}

\def\al{{\alpha}}

\def\om{{\omega}}        
\def\Om{{\Omega}}

\def\si{{\sigma}}

\def\ep{{\varepsilon}}

\def\phi{{\varphi}}

\let\Alpha\alpha
    
\let\pa\partial    
\let\na\nabla    
    
\DeclareMathAlphabet{\doba}{U}{msb}{m}{n}

\gdef\mN{\doba{N}}

\gdef\mR{\doba{R}}        
\gdef\mS{\doba{S}}        
        
\gdef\mZ{\doba{Z}}        
        
\def\lamin{\lambda_{\rm min}}    
  
\def\geucl{g_{\mathrm{eucl}}}    
\def\Vol{{\mathop{\rm Vol}}}

\def\Hom{{\mathop{\rm Hom}}}    
\let\vol\Vol
\def\Scal{{\mathop{\rm Scal}}}
\let\scal\Scal    
    
\def\Supp{{\mathop{\rm Supp}}}    
\def\PSO{{\mathop{P_{\rm SO}}}}    
\def\SO{{\mathop{\rm SO}}}    
    
\def\SU{{\mathop{\rm SU}}}

\let\cM\M
\def\lam{{\lamin^{x_0,y_0}}} 
\def\T{{T_{x_0,y_0}}} 
\let\ti\tilde  
\def\eref#1{{\rm (\ref{#1})}}  

\let\<\langle
\let\>\rangle
\let\embed\hookrightarrow

\long\def\komment#1{}
    
\begin{document}    
%%%%%%%%%%%%%%%%%%%%%%%%%%%%%%%%%%%%%%%%%%%%%%%%%%%%%%%%%%%%%%%%%    
\title{The first conformal Dirac eigenvalue on $2$-dimensional tori}    
\author{Bernd Ammann and Emmanuel Humbert}  
\date{March 2005}  
%\date{\today}    
    
\begin{abstract}    
Let $M$ be a compact manifold with a spin structure $\chi$
and a Riemannian metric $g$. 
Let $\lambda_g^2$ be the smallest eigenvalue of the square of the Dirac operator with respect to $g$ and $\chi$.
The $\tau$-invariant is defined as
  $$\tau(M,\chi):= \sup \inf \sqrt{\lambda_g^2} {\rm Vol}(M,g)^{1/n}$$
where the supremum runs over the set of all conformal classes on $M$, 
and where the infimum runs over all metrics in the given class. 

We  show that $\tau(T^2,\chi)=2\sqrt{\pi}$ 
if $\chi$ is ``the'' non-trivial spin structure on $T^2$.
In order to calculate this invariant, we study the infimum as a function 
on the spin-conformal moduli space and we show that the infimum
converges to $2\sqrt{\pi}$ at one end of the spin-conformal moduli 
space.
\end{abstract}    
  
\maketitle    
\footnote{ammann@iecn.u-nancy.fr, humbert@iecn.u-nancy.fr}   
    
\tableofcontents

\section{Introduction}

Let $(M,g,\chi)$ be a compact spin manifold of dimension $n \geq 2$. 
For any metric $\tilde{g}$ in the conformal class $[g]$ of $g$, 
let $\lambda_1(D_{\tilde{g}}^2)$ be the smallest eigenvalue of the square of 
the Dirac operator. We define     
$$\lamin(M,g,\chi) = \inf_{\tilde{g} \in [g]} \sqrt{\lambda_1(D_{\tilde{g}}^2)}\;   \Vol(M,\tilde{g})^{1/n}. $$    
Several works have been devoted to the study of this conformal   
invariant and some variants of it
\cite{hijazi:86,lott:86,baer:92b,ammann:03,ammann:p03,ammann:habil}.
J.~Lott \cite{lott:86, ammann:03} proved that 
$\lamin(M,g,\chi)=0$ if and only if $\ker D_g\neq \{0\}$.
{}From \cite{hijazi:86,baer:92b} we deduce $\lamin(\mS^n) = \frac{n}{2}\,
\om_n^{\frac{1}{n}}$, where $\mS^n$ is the sphere with constant 
sectional curvature $1$ and where $\om_n$ is its volume. 
Furthermore, in \cite{ammann:03,ammann.humbert.morel:03}, we have seen that
\begin{eqnarray} \label{metric} 
  \lamin(M,g,\chi) \leq   \lamin(\mS^n) = \frac{n}{2}\,
\om_n^{\frac{1}{n}} 
\end{eqnarray} 
for all Riemannian spin manifolds.

Furthermore we define
  $$\tau(M,\chi):=\sup \lamin(M,g,\chi)$$
where the supremum runs over all conformal classes on $M$. Obviously,
$\tau(M,\chi)$ is an invariant of a differentiable manifold with spin 
structure.

We consider it as 
interesting to determine $\tau$ or at least some bounds for $\tau$ 
in as many cases as possible.
There are several motivations for studying these invariants 
$\lamin(M,g,\chi)$ and $\tau(M,\chi)$.

Our first motivation is the analogy and the relation to Schoen's 
$\sigma$-constant, which is defined as
   $$\si(M):=\sup \inf \frac{\int\scal_g \,dv_g}{\Vol(M,g)^{{n-2\over n}}}$$
where the infimum runs over all metrics in a conformal class $g\in [g_0]$,
and where the supremum runs over all conformal classes. 

In the case $\si(M)\geq 0$ and $n\geq 3$, 
there is also an alternative definition
of the $\tau$-invariant that is analogous to our definition of the 
$\tau$-invariant. More exactly, in this case
  $$\si(M):=\sup  \inf \lambda_1(L_g)\Vol(M,g)^{2/n},$$
where  $\lambda_1(L_g)$ is the first eigenvalue
of the conformal Laplacian $L_g:= 4 \,{n-1\over n-2}\,\Delta_g + \Scal_g$.  
Once again, the infimum runs over all metrics in a 
conformal class $g\in [g_0]$,
and where the supremum runs over all conformal classes.
%Schoen's $\sigma$-invariant is defined as the corresponding sup-inf 
%of the first eigenvalue of the conformal Laplacian $L_g:= 4 \,{n-1\over n-2}\,
%\Delta_g + \Scal_g$. 
Many conjectures about the value of the $\si$-constant exist, 
but unfortunately it 
can be calculated only in very few special cases, e.g.\ 
$\si(S^n)={n(n-1)\om_n^{2/n}}$, 
$\si(S^{n-1}\times S^1)={n(n-1)\om_n^{2/n}}$,
$\si(T^n)=0$ and 
$\si(\mR P^3)= {n(n-1)\left({\om_n\over 2}\right)^{2/n}}$. 
The reader might consult \cite{bray.neves:04} for 
a very elegant and amazing calculation of $\si(\mR P^3)$ and for 
a good overview over further literature. 

For other quotients of the sphere $\Gamma\backslash S^n$, 
$\Gamma\subset O(n+1)$
it is conjectured that
\begin{equation}\label{conj.sphere.si}  
  \si(\Gamma\backslash S^n)={n(n-1)\left({\om_n\over \#\Gamma}\right)^{2/n}}.
\end{equation}
It is not difficult to show that for any metric conformal 
to the round metric on $\Gamma\backslash S^n$ one has the inequality
$\lambda_1(L_g)\vol(\Gamma\backslash S^n,g)^{2/n}\geq {n(n-1)\left({\om_n\over \#\Gamma}\right)^{2/n}}.$
This immediately implies $\si(\Gamma\backslash S^n)\geq {n(n-1)
\left({\om_n\over \#\Gamma}\right)^{2/n}}$, i.e.\ the lower bound on $\si$ 
in~\eref{conj.sphere.si}.
However, it is very difficult to obtain the upper bound on $\sigma$.

The $\tau$-invariant is not only a formal analogue to Schoen's 
$\sigma$-constant, but it is also  tightly related to it via
Hijazi's inequality \cite{hijazi:86,hijazi:91,hijazi:01}. Hijazi's inequality  
implies that if $M$ carries a spin structure $\chi$, then 
\begin{equation}\label{ineq.sup.hijazi}
  \tau(M,\chi)^2\geq {n\over 4(n-1)}\,\sigma(M).
\end{equation}
Equality is attained in this inequality if $M=S^n$.
Hence, upper bounds for $\tau(M,\chi)$ may help to determine the 
$\sigma$-constant.

This is one reason for studying the $\tau$-invariant.

%The drawback of this approach is that $\tau(M,\chi)$ is not 
%easier to calculate in general. 

%We know of two other motivations for studying $\lamin(M,g,\chi)$ and 
%the $\tau$-invariant. 
Another motivation for studying $\tau(M,\chi)$ and $\lamin(M,g,\chi)$
comes from the connection to constant mean curvature 
surfaces. Let $n=2$.
If $\tilde g$ is a minimizer that attains the infimum in the definition of
$\lamin(M,g,\chi)$, and if  
$\vol(M,\ti g)=1$, then any simply connected open subset $U$ of $M$
can be isometrically embedded into $\mR^3$, $(U,\ti g)\embed \mR^3$,
such that the resulting surface has constant mean curvature $\lamin(M,g,\chi)$. 
Vice versa, any 
constant mean curvature surface gives rise to a stationary point of an 
associated variational principle.
It is shown in \cite{ammann:p03} that minimizers of $\lamin(M,g,\chi)$ 
exist if $\lamin(M,g,\chi)<2\sqrt{\pi}$.

For the third motivation, let again $n\geq 2$ be arbitrary.
As indicated above, $\lamin(M,g,\chi)>0$ if and only if $\ker D_g=\{0\}$.
Hence, $\tau(M,\chi)>0$ if and only if 
$M$ carries a metric with $\ker D=\{0\}$. It follows from the Atiyah-Singer
index theorem  that any spin manifold $M$
of dimension $4k$, $k\in \mN$ with $\hat A(M)\neq 0$ has $\tau=0$, 
and the same holds for spin manifolds of dimension $8k+1$ and $8k+2$ 
with non-vanishing $\Alpha$-genus.
C. B\"ar conjectures 
\cite{baer:96,baer.dahl:02} that in all remaining cases one has
$\tau>0$. 
%\cite{hitchin:74}, \cite{bourguignon.gauduchon:92}, \cite{baer.schmutz:92},
%\cite{maier:97} . 
Using perturbation methods Maier \cite{maier:97} has verified 
the conjecture in the case $n\leq 4$. 
The conjecture also holds if $n\geq 5$ and
$\pi_1(M)=\{e\}$. Namely, if $M$ is a compact simply connected 
spin manifold with vanishing $\Alpha$-genus, then building on
Gromov-Lawson's surgery results \cite{gromov.lawson:80} Stolz showed
\cite{stolz:92} that $M$ carries a metric $g_+$
of positive scalar curvature. Applying the Schr\"odinger-Lichnerowicz
formula we obtain $\ker D_{g_+}=\{0\}$, and hence 
$\tau(M,\chi)\geq \lamin(M,g_+,\chi)>0$ for the unique spin structure $\chi$
on $M$. A good reference for this argument is also \cite{baer.dahl:02}, 
where the interested reader can also find an analogous statement 
for the case $\Alpha(M,\chi)\neq 0$.
The method of Stolz and B\"ar-Dahl 
also applies to some other fundamental 
groups, but the general case still remains open. 

%It is even more difficult to calculate the actual value of 
%the $\tau$-invariant. 
%Because of this we want to restrict in this
%article to the simplest case, i.e. the case that $M$ is the $2$-dimensional
%torus $T^2$.

In the present article we want to have a closer look at the $\tau$-invariant 
on surfaces, in particular $2$-dimensional tori. 
The higher dimensional case will be the subject of 
another publication.

On surfaces the Yamabe operator cannot be defined as above. 
The Gauss-Bonnet theorem says that the $\si$-constant of a surface does not depend on the metric:
  $$\si(M)=\sup\inf\int 2 K_g\,dv_g= 4\pi \chi(M).$$ 
It was conjectured by Lott \cite{lott:86} and proved by C.~B\"ar \cite{baer:92b} that equation \eref{ineq.sup.hijazi} also holds in dimension $2$. This amounts in showing $\tau(S^2)=2\sqrt{\pi}$. 
If $M$ is a compact orientable surface of higher genus, then
inequality \eref{ineq.sup.hijazi} is trivial. 

%Nevertheless, if the 
%$\alpha$-genus of $(M,\chi)$ vanishes, Lott's inequality and 
%Maier's perturbation result imply $\tau(M,\chi)>0$. 
We will calulate the  
$\tau$-invariant for the $2$-dimensional torus $T^2$. 
The $2$-dimensional torus $T^2$ has $4$ different spin structures. 
The diffeomorphism group ${\rm Diff}(T^2)$ acts on the space of spin structures
by pullback, and the action has two orbits: one orbit consisting of only one
spin structure, the so-called \emph{trivial spin structure} $\chi_{\rm tr}$
and another orbit consisting of three
spin structures. The torus $T^2$ equipped with the trivial spin structure
has non-vanishing $\Alpha$-genus, thus $\tau=0$.
The main result of this article is the following theorem.

\begin{theorem}\label{theo.tau}
Let $\chi$ be a non-trivial spin structure on the $2$-dimensional torus
$T^2$. Then 
  $$\tau(T^2,\chi)=2\sqrt{\pi}\left(=\lamin(\mS^2)\right).$$ 
\end{theorem}

More exactly, for a fixed non-trivial spin structure $\chi$
we will study $\lamin(M,g,\chi)$ 
as a function on the spin-conformal moduli space $\cM$.
We show that it is continuous (Proposition~\ref{prop.cont}), 
and we show that it can be continuously 
extended to the natural $2$-point compactification 
of $\cM$, i.e. the compactification where both ends are compactified 
by one point each.
It will be easy to show that $\lamin(M,g,\chi)\to 0$ 
at one of the ends. 
However, it is much more involved to prove Theorem~\ref{thm.main}
which states that
$\lamin(M,g,\chi)\to  \lamin(\mS^2)=2\sqrt{\pi}$ at the other end. 

It is evident that Theorem~\ref{thm.main} implies Theorem~\ref{theo.tau}.

For the proof of Theorem~\ref{thm.main}, we have to establish a 
qualitative lower bound for the eigenvalues.
One important ingredient in the proof of Theorem~\ref{thm.main} is to study a suitable
covering of the $2$-torus by a cylinder, and to lift a test spinor to this
covering. Using a cut-off argument in a way similar to  
\cite{ammann.baer:02} we obtain a compactly supported test spinor
on the cylinder. After compactifying the cylinder conformally to the sphere
$S^2$, we can use B\"ar's $2$-dimensional version of \eref{ineq.sup.hijazi}, 
to prove $\lamin(M,g,\chi)\to  \lamin(\mS^2)=2\sqrt{\pi}$ at the other end.

Theorem~\ref{thm.main} and Theorem~\ref{theo.tau} should be seen as a 
spinorial analogue of \cite{schoen:91}. In that article, Schoen studies
the Yamabe invariant on the moduli space of
$O(n)$-invariant conformal structures on $S^1\times S^{n-1}$, $n\geq 3$. 
He shows that at one end of this moduli space, the Yamabe invariant converges
to the Yamabe invariant of $\mS^n$, and hence 
$\si(S^1\times S^{n-1})=\si(\mS^n)$. Combining this result with the 
Hijazi inequality and Theorem~\ref{theo.tau}, one obtains
\begin{corollary} Let $n\geq 2$. Then
$$\tau(S^{n-1}\times S^1,\chi)=\left\{
  \begin{matrix}
 0 \hfill & \mbox{if $n=2$ and if $\chi$ is trivial,}\hfill\cr
 {n\over 2}\,\om_n^{1/n}\hfill & \mbox{otherwise.}\hfill
  \end{matrix}
\right.
$$
 \end{corollary}

The structure of the article is as follows.

In section~\ref{moduli}, we define  the spin-conformal moduli $\cM$
space of $2$-tori and recall some well known facts.
In section~\ref{results}, we state and explain our results. In 
sections~\ref{sec.prelim}, we recall some preliminaries which will be useful
for the proof of Theorem~\ref{thm.main}. 
In Section~\ref{sec.proof} the proof is carried out.

{\it Acknowledgement.} The authors want to thank the referee for many useful 
comments.

%%%%%%%%%%%%%%%%%%%%%%%%%%%%%%%%%%%%%%%%%%%%%%%%%%%%%%%%%%%%%%%%%    
\section{The spin-conformal moduli space of $T^2$} \label{moduli}   
%%%%%%%%%%%%%%%%%%%%%%%%%%%%%%%%%%%%%%%%%%%%%%%%%%%%%%%%%%%%%%%%%    
 
At the beginning of this section we will recall the definition of a spin 
structure. We will only give it in the case $n=2$. For more information
and for the case of general dimension we refer to standard text books 
\cite{friedrich:00,lawson.michelsohn:89,roe:88,berline.getzler.vergne:91}. 
More details about the 2-dimensional case can be 
obtained in \cite{ammann.baer:02} and \cite{ammann:diss, baer.schmutz:92}.

Let $(M,g)$ be an oriented surface with a Riemannian metric $g$.
Let $\PSO(M,g)$ denote the set of oriented orthonormal frames over $M$.
The base point map $\PSO(M,g)\to M$ is an $S^1$ principal bundle.
Let $\alpha:S^1\to S^1$ be the non-trivial double covering, i.e.\ $\alpha(z)=z^2$.
A \emph{spin structure on $(M,g)$} is by definition a pair $(P,\chi)$ where
$P$ is an $S^1$ 
principal bundle over $M$ and where $\chi:P\to \PSO(M,g)$ is a double covering,
such that the diagram
\begin{equation}\label{spincompat}
\begin{array}{cccl}
P \times S^1 &\rightarrow& P & \\
& & &\searrow \\
\downarrow\chi\times\alpha & &   \downarrow\chi & \quad M\\
& & &\nearrow \\
\PSO(M) \times S^1 &\rightarrow& \PSO(M) &
\end{array}
\end{equation}
commutes (in this diagram the horizontal flashes denote 
the action of $S^1$ on $P$ and $\PSO(M)$).
By slightly abusing the notation we will sometimes write $\chi$ for the 
spin structure, assuming that the domain $P$ of $\chi$ is implicitly given.
Two spin structures $(P,\chi)$ and $(\ti P,\ti \chi)$ are isomorphic if
there is an $S^1$-equivariant bijection $b:P\to \ti P$ such that $\ti\chi=\chi\circ b$.

If $\ti g=f^2 g$ is a metric conformal to $g$. Then 
$\PSO(M,\ti g)\to \PSO(M,g)$, $(e_1,e_2)\mapsto (fe_2,fe_2)$ defines an 
isomorphism of $S^1$ principal bundles. The pullback of a spin structure
on $(M,g)$ is a spin structure on $(M,\ti g)$.

In a similar way, if $(M_1,g_1)\to (M_2,g_2)$ is an orientation preserving 
conformal map, but not necessarily a diffeomorphism, then
any spin structure on $(M_2,g_2)$ pulls back to a spin structure on  $(M_1,g_1)$.
 
%In our results 
%the sphere $S^2$ and the torus $T^2$, both equipped with arbitrary metrics, 
%will be particularly important. 
\begin{examples}\label{ex.spin}\ \\[-5mm]
\begin{enumerate}[{\rm (1)}]
\item If $g_0$ is the standard metric on $S^2$.
Then $\PSO(S^2,g_0)=\SO(3)$, and the base point map $\SO(3)\to S^2$ is the map 
that associates to a matrix in $\SO(3)$ the first column.
The double cover $\SU(2)\to \SO(3)$ defines a spin structure on $(S,g_0).$ 
\item Let $\ti g$ be an arbitrary metric on $S^2$. After a possible pullback
by a diffeomorphism $S^2\to S^2$ we can write $\ti g= f^2 g_0$.
The pullback of the spin structure given in (1) under the isomorphism
$\PSO(S^2,\ti g)\to \PSO(S^2,g_0)$ defines a spin structure on $(S^2,\ti g)$.
\item Let $g_1$ be a flat metric on the torus $T^2$. Then a parallel frame
gives rise to a (global) section of $\PSO(T^2,g_1)\to T^2$. Hence, this is a 
trivial $S^1$ principal bundle. The trivial fiberwise double covering 
$T^2\times S^1\to T^2\times S^1$, $(p,z)\mapsto (p,z^2)$ defines 
a spin structure on $(T^2,g_1)$, the so-called \emph{trivial spin structure} 
$\chi_{\rm tr}$ on $(T^2,g_1)$
\item If $\ti g$ is an arbitrary metric on $T^2$. Then we can write 
$\ti g=f^2 g_1$ where $g_1$ is a flat metric. As above, the trivial spin 
stucture on $(T^2,g_1)$ defines a spin structure on $(T^2,\ti g)$. 
This spin structure is also called the \emph{trivial spin structure} $\chi_{\rm tr}$.
\item \label{def.cyl}
For $(x_0,y_0)\in \mR^2\setminus\{0\}$  we define
$$Z_{x_0,y_0} =  \mR^2 / \<(x_0,y_0)\>$$ 
where $\<(x_0,y_0)\>$ is the subgroup of $\mR^2$ spanned by 
$(x_0,y_0)$. We will assume that it carries the metric induced by the 
euclidean metric $\geucl$ on $\mR^2$. Then $\PSO(Z_{x_0,y_0})$ is a trivial bundle, and a natural
trivialization is obtained by a parallel frame.  The map
$Z_{x_0,y_0}\times S^1\to Z_{x_0,y_0}\times S^1$, $(p,z)\mapsto (p,z^2)$ defines
a spin structure on $Z_{x_0,y_0}$, the \emph{trivial spin structure} on $Z_{x_0,y_0}$.
\end{enumerate}
\end{examples}

Assume that $\chi:P\to \PSO(M,g)$ is a spin structure on a surface, and assume
that $\beta:\pi_1(M)\to \{-1,+1\}$ is a group homomorphism. Then there
is a $ \{-1,+1\}$ principal bundle $B_\beta\to M$ with holonomy $\beta$.
Let $P_\beta$ be the quotient of $P\times B_\beta$ by the diagonal action
of $\{-1,+1\}$. Then $P_\beta$ together with the induced  map $\chi_\beta:P_\beta\to \PSO(M,g)$
is also a spin structure on $(M,g)$.  
Conversely, if  $(\ti P,\ti \chi)$ is another spin structure, 
then one can show that there is a unique $\beta:\pi_1(M)\to \{-1,+1\}$ 
such that $(\ti P,\ti \chi)$ and $(P_\beta,\chi_\beta)$ are isomorphic. 
Thus, we see that the space of spin structures is an affine space
over the $\{-1,+1\}$-vector space $\Hom(\pi_1(M),\{-1,+1\})=H^1(M,\{-1,+1\})$.

\begin{examples}\ \\[-5mm]
\begin{enumerate}[{\rm (1)}]
%\item The only homomorphism $\pi_1(S^2)\to \{-1,+1\}$ is the trivial map. Hence, there is only one spin
%structure on $S^2$ up to isomorphism.
\item Any compact oriented surface $M$ carries a spin structure.
If $k$ denotes the genus of $M$, then there are 
$4^k$ homomorphisms $\pi_1(M)\to\{-1,+1\}$,
hence there are $4^k$ isomorphism classes of spin structures. 
In particular, the spin structure
on $S^2$ is unique.
\item Because of $\pi_1(Z_{x_0,y_0})=\mZ$, there are exactly two spin structures on $Z_{x_0,y_0}$, the trivial
one and another one called the \emph{non-trivial spin structure}.
\end{enumerate}
\end{examples}

{}From now  on, let $M=T^2=\mR^2/\Gamma$ where $\Gamma$ is a lattice in $\mR^2$. 
The trivial spin structure defined above can be used to identify 
$\Hom(\pi_1(M),\{-1,+1\})$ with the set of isomorphism classes of spin structures. 
By slightly abusing the language we will always write $\chi$ for the spin structure $(P,\chi)$ and also for the 
homomorphism $\pi_1(M)\to\{-1,+1\}$.

The following lemma summarizes
some well-known equivalent characterizations of triviality of $\chi$
(see e.g.\ \cite{lawson.michelsohn:89}, \cite{milnor:63}, \cite{ammann:diss}, \cite{friedrich:84}).

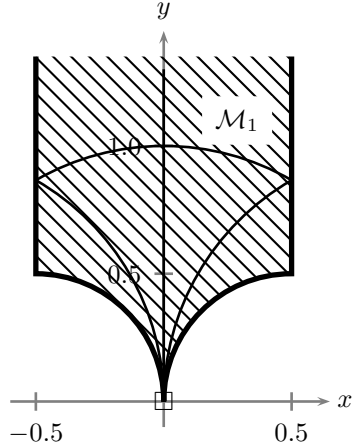
\begin{figure}[tb]
%
% Dies ist die Datei spinmodu.pic. 
% Dieses Bild zeigt den spin-konformen Modulraum
% Teil des Artikel "Spectral estimates on 2-tori"
% New York, 8.12.1999
% 
%
\begin{center}
\newgray{abshilflinfarbe}{.30}
%\newrgbcolor{absliyaufarbe}{1 0 0}  %bunt
%\newrgbcolor{absammannfarbe}{0 1 0} %bunt
\newgray{absliyaufarbe}{.80}  %schwarzweiss
\newgray{absammannfarbe}{.40}  %schwarzweiss

\newdimen\descwidth
\descwidth=6cm
\def\achseneinst{\psset{linewidth=1pt,linecolor=gray,linestyle=solid,fillstyle=none}}

\def\abspinpfeileinst{\psset{linewidth=1pt,linecolor=black,linestyle=solid,fillstyle=none}}
\def\absmodformeinst{\psset{linewidth=2pt,linecolor=black,linestyle=solid,fillstyle=vlines,hatchcolor=black}}
\def\absliyaueinst{\psset{linewidth=2pt,linecolor=black,linestyle=solid,fillstyle=vlines*,fillcolor=absliyaufarbe,hatchcolor=black}}
\def\absammanneinst{\psset{linewidth=2pt,linecolor=black,linestyle=solid,fillstyle=vlines*,fillcolor=absammannfarbe,hatchcolor=black}}
\def\absralilineinst{\psset{linewidth=1pt,linecolor=black,linestyle=solid,fillstyle=none}}
\def\cliffordeinst{\psset{linewidth=2pt,linecolor=black,linestyle=solid,fillstyle=solid,fillcolor=black}}

\psset{unit=1.7cm}
\begin{pspicture}(-1.5,-.2)(1.7,3.35)
\achseneinst
\psaxes[Dx=.5,dx=1,Dy=.5,dy=1]{->}(0,0)(-1.2,-.2)(1.3,2.9)

\absmodformeinst
\pscustom{
  \psline(-1,2.7)(-1,1)
  \psarcn(-1,0){1}{90}{0}
% \psline(0,2.7)(0,0)
  \psarcn(1,0){1}{180}{90}
  \psline(1,1)(1,2.7)
} % of \pscustom

%\absliyaueinst
%\pscustom{
%  \psarc(0,1){1}{0}{90}
%  \psline(0,2)(1,2)
%  \psline(1,2)(1,1)
%}% of \pscustom

%\absammanneinst
%\pscustom{
%  \psline(0,1)(0,0)
%  \psarcn(1,0){1}{180}{90}
%  \psline(1,1)(0,1)
%} % of \pscustom

\absralilineinst
\psarc(0,0){2}{60}{120}
\psarc(2,0){2}{120}{180}
\psarc(-2,0){2}{0}{60}
\psline(0,0)(0,2.6)

\rput(0,0){$\square$}

\cliffordeinst
%\pscircle(1,1){.06}

\rput(1.42,0){$x$}
\rput(0,3.06){$y$}

% \rput[l](1.2,1){\vtop{\hsize=3cm\noindent\raggedright A quadratic\\torus}}
\rput[l](.3,2.2){\psframebox[linecolor=white,fillstyle=solid,fillcolor=white]{$\cM_1$}}

%\def\tere{2.0}
%\def\terz{2.4}
%\rput[r](\tere,1.5){\psframe[linecolor=absliyaufarbe,fillstyle=solid,fillcolor=absliyaufarbe](0,-.15)(.3,.15)}
%\rput[l](\terz,1.5){ 
%\vtop{\hsize=\descwidth\noindent\raggedright Li and Yau proved the Willmore conjecture for these spin-conformal classes}
%}
%\rput[r](\tere,.5){\psframe[linecolor=absammannfarbe,fillstyle=solid,fillcolor=absammannfarbe](0,-.15)(.3,.15)}
%\rput[l](\terz,.5){
%\vtop{\hsize=\descwidth\noindent\raggedright For these spin-conformal classes we prove the Willmore conjecture under a curvature assumption} 
%}
%\absmodformeinst
%\rput[r](2.0,2.5){\psframe(0,-.15)(.3,.15)}
%\rput[l](2.4,2.5){
%\vtop{\hsize=\descwidth\noindent\raggedright
%The spin conformal moduli space $\spinmod$}
%}
\end{pspicture}
\end{center}
\caption{The spin conformal moduli space is $\cM=\cM_1/\sim$, where $\sim$ 
means identifying $(x,y)\in\pa\cM_1$ with $(-x,y)$.}\label{spinmodbild}
\end{figure}

\begin{lemma}
With the above notations, the following statements are equivalent
\begin{enumerate}[\rm (1)]
\item The spin structure is trivial (in the above sense);
\item $\chi(\gamma)=1$ for all $\gamma\in \Gamma$;
\item The spin structure is invariant under the natural action
of the diffeomorphism group ${\rm Diff}(T^2)$;
\item $(T^2,\chi)$ is the {\bf non}-trivial element in the $2$-dimensional 
spin-cobordism group;
\item The $\Alpha$-genus of $(T^2,\chi)$ is the non-trivial element in $\mZ/2\mZ$;
\item The Dirac operator has a non-trivial kernel;
\item The kernel of the Dirac operator has complex dimension $2$.
\end{enumerate}
\end{lemma}

%The Dirac operator with respect to the trivial spin structure and with respect 
%to any metric has a kernel of complex dimension $2$. 
%Hence, 
In particular, we easily see 
  $$\tau(T^2,\chi_{\rm tr})=0.$$

{}From now on, in the rest of this article, we assume that $\chi$ is not the 
trivial spin structure,
%Our article will study the first positive eigenvalue of the Dirac operator
%on tori with a non-trivial spin structure,
i.e.\ $\chi(\gamma)=-1$ 
for some $\gamma\in \Gamma$.

\begin{definition}
Two $2$-dimensional tori with Riemannian metrics, orientations and spin structures
are said to be \emph{spin-conformal} if there is a conformal map between them
preserving the orientation and the spin structure. 
``Being spin-conformal'' is obviously an equivalence relation.
The \emph{spin-conformal moduli space $\cM$ of $T^2$ 
with the non-trivial spin structure} is defined to be the set of 
these equivalence classes.
Furthermore we define
  $$\cM_1:=\left\{
\begin{pmatrix}x\cr y\end{pmatrix}\,\Big|\, |x|\leq {1\over 2},\quad 
y^2+\left(|x|-{1\over 2}\right)^2\geq {1\over 4},\quad y> 0\right\}$$
(see also Fig.~\ref{spinmodbild}).
\end{definition}

\begin{lemma}
Let $g$ be a Riemannian metric on $T^2=\mR^2/\mZ^2$, and let 
$\chi:\mZ^2\to \{-1,+1\}$ be a non-trivial spin structure.
Then there is a lattice $\Gamma\subset \mR^2$,
a spin structure $\chi':\Gamma\to\{-1,+1\}$,
such that 
\begin{enumerate}[{\rm (1)}]
\item $\Gamma$ is generated by
$\begin{pmatrix}1\cr 0\end{pmatrix}$ and 
$\begin{pmatrix}x\cr y\end{pmatrix}$ with $\begin{pmatrix}x\cr y\end{pmatrix}\in \cM_1$
\item $(T^2,g,\chi)$ is spin-conformal to $(\mR^2/\Gamma,g_{\rm eucl},\chi')$
\item $\chi'(\begin{pmatrix}1\cr 0\end{pmatrix})=+1$ and $\chi'(\begin{pmatrix}x\cr y\end{pmatrix})=-1.$
\end{enumerate}
\end{lemma}

\komment{
\begin{lemma}
Let $g$ be a Riemannian metric on $T^2=\mR^2/\mZ^2$, and let 
$\chi:\mZ^2\to \{-1,+1\}$ be a non-trivial spin structure.

Then there is a lattice $\Gamma\subset \mR^2$,
a spin structure $\chi':\Gamma\to\{-1,+1\}$, and an orientation preserving
conformal map 
$F:(T^2,g)\to \mR^2/\Gamma,\geucl$
such that 
\begin{enumerate}[{\rm (1)}]
\item $\chi=\chi'\circ F$
\item $\Gamma$ is generated by
$\begin{pmatrix}1\cr 0\end{pmatrix}$ and 
$\begin{pmatrix}x\cr y\end{pmatrix}$ with $\begin{pmatrix}x\cr y\end{pmatrix}\in \cM_1$
\item $\chi'(\begin{pmatrix}1\cr 0\end{pmatrix})=+1$ and $\chi'(\begin{pmatrix}x\cr y\end{pmatrix})=-1.$
\end{enumerate}
\end{lemma}
}

\proof{}
Because of the uniformization theorem we can assume without loss 
of generality that $g$ is a flat metric.  The lemma then follows from 
elementary algebraic arguments.
\qed

Note that $x$ and $y$ are uniquely determined if 
$\begin{pmatrix}x\cr y\end{pmatrix}$ is in the interior of $\cM_1$, i.e. if
$|x|< 1/2$ and $y^2+(|x|-1/2)^2> 1/4$. 
If $\begin{pmatrix}x\cr y\end{pmatrix}$ is on the boundary of $\cM_1$, 
then $y$ and $|x|$ are determined, 
but not the sign of $x$. Hence, after gluing 
$\begin{pmatrix}x\cr y\end{pmatrix}\in \pa \cM_1$ with 
$\begin{pmatrix}-x\cr y\end{pmatrix}$ we obtain the 
spin-conformal moduli space $\cM$.

{\it Notation.}  
Let $(x_0,y_0) \in \cM_1$. The lattice generated by 
 by
$\begin{pmatrix}1\cr 0\end{pmatrix}$ and 
$\begin{pmatrix}x\cr y\end{pmatrix}$
is noted as $\Gamma_{x_0,y_0}$. Furthermore, we write 
$T_{x_0,y_0}$ for the $2$-dimensional torus $\mR^2/\Gamma_{x_0,y_0}$
\emph{equipped with the euclidean metric.}

The quantity $\lamin(T^2,g,\si)$ is a spin-conformal invariant, hence $\lamin$
can be viewed as a function on $\cM$ or on $\cM_1$.

\section{Main results} \label{results}

%In this article, 
%we are interested in the function  
%\[ \lamin: \left| \begin{array}{ccc} 
%\cM_1 & \to &  ]0, \lamin(\mS^2)] \\[2mm] 
%(x_0,y_0) & \mapsto & \lam. 
%\end{array} \right. \]
In this article, we  study $\lamin$ as a function
on the spin-conformal moduli space with the non-trivial spin structure. 
This function takes values in $[0,\lamin(\mS^2)]$ because of \eref{metric}.
As the spin structure is non-trivial, Lott's results states that $0$ 
is not attained.
As a preliminary result we will prove
that this function is continuous.

\begin{proposition} \label{prop.cont} 
The function  
\[ \lamin: \left| \begin{array}{ccc} 
\cM_1 & \to &  ]0, \lamin(\mS^2)] \\[2mm]
(x_0,y_0) & \mapsto & \lam
\end{array} \right. \] 
is continuous on $\cM_1$. 
\end{proposition} 
 
The spin-conformal moduli space $\cM$ (resp. $\cM_1$) has two ends.
We will compactify each end by adding one point.
The point added at the end $y\to \infty$ will be denoted by $\infty$
and the point added at the end $y\to 0$ is denoted by $(0,0)$.
 
\begin{thm} \label{thm.main} 
The function  
\[ \lamin: \left| \begin{array}{ccc} 
\cM_1 & \to &  ]0, \lamin(\mS^2)] \\ 
(x_0,y_0) & \mapsto & \lam
\end{array} \right. \] 
extends continuously to  
$\cM_1 \cup \{(0,0),\infty\}$ by setting $ \lamin^{0,0} =  \lamin(\mS^2)$
and $\lamin^\infty=0$. 
\end{thm} 
 
The continuous extension at $\infty$ is is easy to see.
The first eigenvalue of the Dirac operator on $(T_{x_0,y_0},g_{eucl},
  \chi_{x_0,y_0})$, is $\pi/y_0$, the area is $y_0$, hence 
  $$ \lam\leq \pi/\sqrt{y_0}\to 0\qquad \mbox{for }y_0\to \infty.$$ 

However, the limit $(x_0,y_0)\to(0,0)$ is much more 
difficult to obtain.
 
Clearly, Theorem~\ref{thm.main} implies Theorem~\ref{theo.tau}.

\section{Some preliminaries}\label{sec.prelim}

{\bf Variational characterization of $\lamin$.} %\label{functional}
Let $(M,g,\chi)$ be a compact spin manifold of dimension $n \geq 2$
with $\ker D_g=\{0\}$. For 
$\psi\in\Gamma(\Sigma M)$, we define     
$$J_g(\psi)=\frac{\Big(\int_M|D\psi|^{\frac{2n}{n+1}}\,
dv_g\Big)^\frac{n+1}{n}}{|\int_M
  \<D\psi,\psi\>\,dv_g|} .$$  
Lott \cite{lott:86} proved  that    
\begin{eqnarray} \label{func}   
\lamin(M,[g],\chi)=\inf_\psi J_g(\psi)   
\end{eqnarray}   
where the infimum is taken over the set of smooth spinor fields for which    
$$\left(\int_M  \<D\psi,\psi\>\,dv_g \right)\neq 0.$$

The functional $J_g$ for the torus $T_{x_0,y_0}$ is noted as $J^{x_0,y_0}$. 

\begin{remark}
The exponents in $J_g$ are chosen such that  
$J_g$ is conformally invariant. More exactly, if $g$ and $\ti g$ are conformal,
then the spinor bundles of $(M,g,\chi)$ and $(M,\ti g, \chi)$ can be 
identified in such a way that $J_g(\psi)=J_{\ti g}(\psi)$.
\end{remark}

{\bf Cylinders and doubly pointed spheres.} 

Let $Z_{x_0,y_0}$ be defined as in Examples~\ref{ex.spin} (\ref{def.cyl}).

\begin{lem}[Mercator, around 1569]\label{lem1}
Let $N,S \in \mS^2$ be respectively the North pole and the South pole of
$\mS^2$. Then there is a conformal diffeomorphism $F_{x_0,y_0}$ from 
$(Z_{x_0,y_0}, g_{eucl})$ to $(\mS^2 \setminus \{N,S\})$.
\end{lem} 

\proof{}
In the case $(x_0,y_0)=(0,2\pi)$ we see that
the application
$$F_{0,2\pi}:\begin{pmatrix}x\cr y\end{pmatrix}\mapsto 
  \begin{pmatrix}{\sin y \over \cosh x}\cr{\cos y\over \cosh x}\cr \tanh x\end{pmatrix}$$
is conformal and defines a conformal bijection 
$Z_{0,2\pi} \to S^2\setminus\{N,S\}$. 
The general case follows by composing with a 
linear conformal map $Z_{x_0,y_0}\to Z_{0,2\pi}$.
\qed

The map $F$ induces a map between the frame bundles.
\begin{align*}
  \ti F_{x_0,y_0}&:\PSO(Z_{x_0,y_0})\to \PSO(\mS^2)\\ 
  \ti F_{x_0,y_0}((p,X,Y))&:=\left(F_{x_0,y_0}(p),{dF_{x_0,y_0}(X)\over |dF_{x_0,y_0}(X)|}, {dF_{x_0,y_0}(Y)\over |dF_{x_0,y_0}(Y)}\right)\\
   X,Y&\in T_p Z_{x_0,y_0} \mbox{ are orthonormal and oriented}
\end{align*}
The unique spin structure on $\mS^2$ pulls back to a spin structure
on $Z_{x_0,y_0}$, that we will denote as $\chi_{x_0,y_0}$.

\begin{lemma}
The spin structure $\chi_{x_0,y_0}$ is the non-trivial spin structure 
on $Z_{x_0,y_0}$.
\end{lemma}

\proof{} We will show the lemma for the case $(x_0,y_0)=(0,2\pi)$. 
As before, the general case then follows by composing with a linear map
 $Z_{x_0,y_0}\to Z_{0,2\pi}$.

We define the loop 
$\gamma:[0,2\pi]\to Z_{0,2\pi}$, $\gamma(t):=(0,t)$ and the parallel section 
$$\al:t\mapsto ({\pa\over\pa x}|_{\gamma(t)},{\pa\over\pa y}|_{\gamma(t)})$$
of $\PSO(Z_{0,2\pi})$ along $\gamma$.
The spin structure $(P,\chi_{0,2\pi})$
on  $Z_{0,2\pi}$ is  trivial if and only if there
is a section  $\ti\al$ of $P$ along $\gamma$ such that 
$\chi_{0,2\pi}\circ \ti\al=\al$ and $\ti\al(0)=\ti\al(2\pi)$.

The composition $\ti F_{0,2\pi}\circ \al$ is a section of 
$\PSO(\mS^2)=\SO(3)$ along $F_{0,2\pi}\circ \gamma$. 
One checks that  
  $$\ti F_{0,2\pi}\circ \al(t)
    =\left({\pa F_{0,2\pi}\over \pa x}|_{(0,t)},
           {\pa F_{0,2\pi}\over \pa y}|_{(0,t)}, F(0,t)\right)=
   \begin{pmatrix}
      0 & \cos y & \sin y \\
      0 & -\sin y & \cos y \\
      1 & 0 & 0
   \end{pmatrix}
  $$
We lift this loop to a path $\hat\al$ in $\SU(2)$, then one easily sees
that $\hat\al(0)=-\hat\al(2\pi)$. As  
$\chi_{0,2\pi}$ is defined as the pullback of the spin structure on $S^2$,
we see that any lift $\ti\al$ of $\al$ also satisfies
$\ti\al(0)\neq \ti\al(2\pi)$.
Hence, we have proved non-triviality of $\chi_{0,2\pi}$.
\qed

\begin{cor} \label{cor2}
Let $Z_{x_0,y_0}$ carry its non-trivial spin structure. Then, 

$$ \frac{{\left( \int_{Z_{x_0,y_0}} {|D \psi|}^{\frac{4}{3}} dx
    \right)}^{\frac{3}{2}} }{ \left| \int_{Z_{x_0,y_0}} \big\< \psi, D \psi\big\>
  dx\right|} \geq \lamin(\mS^2)$$
for any compactly supported
spinor $\psi \in \Gamma( \Sigma Z_{x_0,y_0})$ such that 
$\int\<\psi,D\psi\>\neq 0$.
\end{cor} 
%Indeed, by Lemma~\ref{lem1}, there exists a conformal diffeomorphism 
%$\theta:(Z_{x_0,y_0}, g_{eucl}) \to (\mS^2 \setminus
%\{N,S\},g_0)$. 
Let $f: Z_{x_0,y_0} \to ]0, +\infty[$ be such that 
$F_{x_0,y_0}^* g_0 = f^2 g_{eucl}$. 
It is well known 
(see for example \cite{hitchin:74,hijazi:86} ) that 
$F_{x_0,y_0}$ induces a pointwise isometry 
\[ \left| \begin{array}{ccc}
\Sigma( \T, g_{eucl} ) & \to & \Sigma (\mS^2 \setminus 
\{N,S\},g_0) \\  
 \psi & \mapsto&  \bar{\psi},
\end{array} \right. \]
such that 
  $$\bar{D}  f^{-\frac{1}{2}} \bar{\psi} 
     = f^{-\frac{3}{2}} \overline{D \psi}.$$
%Let now $\psi \in \Gamma( \Sigma Z_{x_0,y_0})$ be compactly supported in
%$Z_{x_0,y_0}$. 
%Moreover, if we set 
%$$\phi= f^{-\frac{1}{2}} \bar{\psi}$$
where $\bar{D}$ denotes the Dirac operator on $\mS^2$.
Moreover, $\bar{\psi}$ is smooth on $\mS^2$ since
$\psi \equiv 0$ in a neighborhood of $N$ and $S$. It is well known that
the functional $J$ defined  
at the beginning of section~\ref{sec.prelim} is 
conformally invariant. This implies that
$$ \frac{{\left( \int_{Z_{x_0,y_0}} {|D \psi|}^{\frac{4}{3}}\, dx
    \right)}^{\frac{3}{2}} }{ | \int_{Z_{x_0,y_0}} \big\< \psi, D \psi\big\>\,
  dx|}  =  \frac{{\left( \int_{\mS^2} {|\bar{D}  (f^{-\frac{1}{2}} \bar{\psi})|}^{\frac{4}{3}}\, dv_{g_0}
    \right)}^{\frac{3}{2}} }{ | \int_{\mS^2} 
   \big\<  f^{-\frac{1}{2}} \bar{\psi}, \bar{D} 
    (f^{-\frac{1}{2}} \bar{\psi})\big\>\,
    dv_{g_0}|} \geq \lamin(\mS^2).$$

\section{Proof of the Main results}\label{sec.proof}

For the proof we will need the following well known elliptic estimates.
These estimates are a consequence of techniques explained
for example in \cite{taylor:81}, see also
\cite{aubin:98}.
However, in our special situation a proof is much easier. Hence, for the 
convenience of the reader we will include an elementary proof here. 

\begin{lemma}[Elliptic estimates]\label{lem.ell}
Let $(x_0,y_0) \in \cM_1$, and note $T^2$ for $T_{x_0,y_0}$. 
%We use the coordinates $(T^2,g=g_{x_0,y_0},\chi)$ (see paragraph \ref{view}). 
There exists $C>0$
depending only on $x_0$ and $y_0$ such that 
\begin{eqnarray} \label{r4}
\int_{T^2} |D \psi |^{\frac{4}{3}}  dv_{g} \geq C \int_{T^2} | \nabla \psi
|^{\frac{4}{3}}   dv_{g}
\end{eqnarray}
and 
\begin{eqnarray} \label{r5}
{\left( \int_{T^2} |\psi|^{4} dv_g\right)}^{\frac{1}{3}}
\leq   C \int_{T^2} | \nabla \psi
|^{\frac{4}{3}}   dv_{g}
\end{eqnarray}
for any smooth spinor $\psi$.
\end{lemma}

\proof{}
Let $q=\frac{4}{3}$.
Assume that (\ref{r4}) is false. Then,  for all $\ep >0$,  we can find a 
 smooth
spinor $\psi_\ep \in \Gamma(\Sigma(T^2))$ such 
that  
\begin{eqnarray} \label{boun}
\int_{T^2} |D \psi_\ep |^q  \,dv_{g} \leq \ep \qquad \mbox{and} \qquad
\int_{T^2} | \nabla \psi_\ep |^q \, dv_{g}=1.
\end{eqnarray} 
Now, assume that  
$$\lim_{\ep \to 0} {\left( \int_T |\psi_\ep|^q\, dv_g \right)}^{\frac{1}{q}}=
+\infty.$$
Then, we set  
$$ \psi'_\ep= \frac{\psi_\ep}{ {\left( \int_{T^2} |\psi_\ep|^q dv_g
    \right)}^{\frac{1}{q}}}.$$
The sequence $(\psi'_\ep)$ is bounded in $W^{1,q}(T^2)$ and since
$W^{1,q}(T^2)$ is reflexive, 
we can find $\psi'_0  \in W^{1,q}(T^2)$ such that there is sequence $\ep_i\to 0$,
with $\lim_{i\to \infty} \psi'_{\ep_i}=\psi'_0 $ weakly in $W^{1,q}(T^2)$. 
Then, we would have
$$ \int_{T^2}  |\nabla \psi'_0 |^q   dv_g \leq \liminf_{\ep}  \int_{T^2}  |\nabla
\psi'_\ep |^q   dv_g =0$$
We would get  that $\psi'_0$ is parallel which cannot occur since the
structure on $T^2$ is not trivial.
This proves that $(\psi_\ep)$ is bounded in $L^q(T^2)$ and hence,
by (\ref{boun}) in
$W^{1,q}(T^2)$. 
Again by reflexivity of $W^{1,q}(T^2)$, we get  the existence of a spinor 
$\psi_0 $, weak limit of a subsequence $\psi_{\ep_i}$ in $W^{1,q}(T^2)$. 
By weak convergence
of $D\psi_{\ep_i} $ to $D \psi_0$ in $L^q(T^2)$, we have  
$$\int_{T^2} |D \psi_0 |^q dv_g \leq \liminf_i  \int_{T^2} |D \psi_{\ep_i}  |^q dv_g
= 0.$$
This is impossible since the Dirac operator on $T^2$ has a trivial
kernel. 
This proves (\ref{r4}). As one can check, relation (\ref{r5}) can be proved
with the same
type of arguments. 
\qed

\proof{of Proposition~\ref{prop.cont}}
The proposition states that  $\lamin$ is continuous on $\cM_1$. 
Let $(x_k,y_k)_k\in \cM_1$ be a sequence tending to $(x_0,y_0)\in \cM_1$.
We identify $T^2$ with $\mR^2/\mZ^2$.
The conformal structures corresponding to $(x_k,y_k)$ and $(x_0,y_0)$
are represented by flat metrics  $g_{x_k,y_k}$ and $g_{x_0,y_0}$
on  $\mR^2/\mZ^2$, that are invariant under translations, and such that
$g_{x_k,y_k}\to g_{x_0,y_0}$ in the $C^\infty$-topology . 

Let $\ep >0$ be small and  let $\psi_0$ and 
$(\psi_k)_k$ be smooth spinors such that 
$$J_{x_0,y_0} (\psi_0) \leq \lam + \ep \hbox{ and } J_{x_k,y_k} (\psi_k) \leq
\lamin^{x_k,y_k}+ \ep. $$
At first, since $(g_{x_k,y_k})_k$ tends to $g_{x_0,y_0}$, it is easy to
see that 
$$\lim_k J_{x_k,y_k} (\psi_0) = J_{x_0,y_0} (\psi_0)$$
and hence 
$\limsup_k \lamin^{x_k,y_k} \leq \lam+ \ep$
for the given $\ep>0$ that we can choose as small as we want. Thus
$$\limsup_k \lamin^{x_k,y_k} \leq \lam.$$
Now, let us prove that

\begin{eqnarray} \label{r2}
\limsup_k J_{x_0,y_0} (\psi_k) \leq \liminf_k  J_{x_k,y_k} (\psi_k)
%\leq \liminf_k \lamin^{x_k,y_k} +\ep.
\end{eqnarray}

\noindent We let $(v,w)$ be a orthormal basis for $g_{x_0,y_0}$ and 
$(v_k,w_k)_k$,  orthonormal basis for $g_{x_k,y_k}$ which tends to
$(v,w)$. One can write for all
$k$, $v_k=a_k v+b_k w$ and $w_k= c_k v + d_k w$ with 
$\lim_k a_k= \lim_k d_k=1$ and $\lim_k b_k= \lim_k c_k=0$. We have 
$${\left(\int_{T^2} |D_{x_k,y_k}\psi_k |^{\frac{4}{3}} 
dv_{g_{x_k,y_k}} \right)}^{\frac{3}{4}} = 
{\left(\int_{T^2} |v_k \cdot \nabla_{v_k}\psi_k + w_k \cdot \nabla_{w_k} \psi_k |^{\frac{4}{3}} 
dv_{g_{x_k,y_k}} \right)}^{\frac{3}{4}}$$
$$= {\left(\int_{T^2} | D_{x_0,y_0} \psi_k + \theta_k
    |^{\frac{4}{3}}  
dv_{g_{x_k,y_k}} \right)}^{\frac{3}{4}}.$$
with 
$$|\theta_k| = |(a_k^2+c_k^2- 1 ) v \cdot \nabla_v \psi_k + 
(b_k^2+d_k^2- 1 )  w \cdot \nabla_w \psi_k
+ (a_k b_k+c_k d_k ) \<w \cdot \nabla_v \psi_k +  v \cdot \nabla_w \psi_k\>
| \leq \alpha_k |\nabla \psi_k| $$
where $(\alpha_k)_k$ is a sequence of positive numbers which tends to
$0$. Note that because of the translation invariance of the metrics,
the Levi-Civita connection does not depenpend on $k$.
Since $\lim_k g_{x_k,y_x} =g_{x_0,y_0}$, one gets that 
$${\left(\int_{T^2} |D_{x_k,y_k}\psi_k |^{\frac{4}{3}} 
dv_{g_{x_k,y_k}} \right)}^{\frac{3}{4}} 
\geq {\left(\int_{T^2} |D_{x_0,y_0}\psi_k |^{\frac{4}{3}} 
dv_{g_{x_0,y_0}} \right)}^{\frac{3}{4}} - \alpha_k' 
{\left( \int_{T^2}  |\nabla \psi_k|^{\frac{4}{3}} 
dv_{g_{x_0,y_0}} \right)}^{\frac{3}{4}}$$
where $\lim_k  \alpha'_k =0$.
Together with Lemma~\ref{lem.ell}, we get that 
\begin{eqnarray} \label{r3}
(1- C \alpha_k') {\left(\int_{T^2} |D_{x_0,y_0}\psi_k |^{\frac{4}{3}} 
dv_{g_{x_0,y_0}} \right)}^{\frac{3}{4}} \leq 
 {\left(\int_{T^2} |D_{x_k,y_k}\psi_k |^{\frac{4}{3}} 
dv_{g_{x_k,y_k}} \right)}^{\frac{3}{4}}  
\end{eqnarray}
where $C$ is a positive constant independent of $k$. 
Now, in the same way, we can write 
 
$$\int_{T^2}  \big\<\psi_k, D_{x_0,y_0} \psi_k \big\>\, dv_{g_{x_0,y_0}} \geq 
\int_{T^2}  \big\<\psi_k, D_{x_k,y_k} \psi_k \big\>\, dv_{g_{x_k,y_k}}
- \beta_k  \int_{T^2} |\psi_k|  |\nabla \psi_k| dv_{g_{x_0,y_0}}$$
where $\lim_k \beta_k= 0 $. Using H\"older inequality, we have 
$$\int_{T^2} |\psi_k| |\nabla \psi_k| dv_{g_{x_0,y_0}}
\leq {\left( \int_{T^2}  |\psi_k|^{4} dv_{g_{x_0,y_0}} \right)}^{
{\frac{1}{4}}}  {\left( \int_{T^2}  |\nabla  \psi_k|^{\frac{4}{3}}
    dv_{g_{x_0,y_0}} \right)}^{ 
{\frac{3}{4}}}. $$
Using (\ref{r4}) and (\ref{r5}), this gives 
$$\int_{T^2} |\psi_k| |\nabla \psi_k| dv_{g_{x_0,y_0}} 
\leq C  \left(\int_{T^2}  |D  \psi_k|^{\frac{4}{3}}
    dv_{g_{x_0,y_0}}\right)^{3/2}.$$ 
We obtain 

$$\int_{T^2}  \big\<\psi_k, D_{x_0,y_0} \psi_k \big\>\, dv_{g_{x_0,y_0}} \geq 
\int_{T^2}  \big\<\psi_k, D_{x_k,y_k} \psi_k \big\>\, dv_{g_{x_k,y_k}}
- \beta_k  \left(\int_{T^2} |D \psi_k|^{\frac{4}{3}}\, dv_{g_{x_0,y_0}}.
 \right)^{3/2}$$
Together with (\ref{r3}), we get (\ref{r2}). This immediatly implies 
that 
$$\liminf_k \lamin^{x_k,y_k} \geq \lam$$
and ends the proof of the proposition.\qed

\proof{of Theorem~\ref{thm.main}}
Any calculation in this proof will be carried out in Riemannian normal
coordinates with respect to a flat metric. In the following, $(e_1,e_2)$ 
will denote the canonical basis of $\mR^2$.

In order to prove 
$\lim_{(x_0,y_0) \to (0,0)} \lam = \lamin(\mS^2)$
we will show that there is no sequence $(x_k,y_k)\to (0,0)$ such that 
$\lim_{(x_k,y_k) \to (0,0)}\lamin^{x_k,y_k} < \lamin(\mS^2)$.
We may assume that $\lamin^{x_k,y_k} < \lamin(\mS^2)$ for all $k$.
%We use the coordinates $(\T,g_{eucl})$.
Note that the spectrum of $D$ is symmetric in dimension $2$.
By \cite{ammann:habil}, we then can find a
sequence of spinors $\psi_k$ of class $C^1$ such that on $T_{x_k,y_k}$

\begin{eqnarray} \label{eqk}
D \psi_k = \lamin^{x_k,y_k} |\psi_k|^2 \psi_k
\end{eqnarray}
and such that
\begin{eqnarray} \label{norm}
\int_{T_{x_k,y_k}} |\psi_k|^{4} dx = 1.
\end{eqnarray}
Moreover, we have 
\begin{eqnarray} \label{opt}
J_{x_k,y_k}( \psi_k) = \lamin^{x_k,y_k}.
\end{eqnarray}
Sometimes we will identify $\psi_k$ with its pullback to $\mR^2$.
In this picture $\psi_k$ is a doubly periodic spinor on $\mR^2$.

\begin{step} \label{st2}  
There exists $C>0$ such that for all   $k$, we have   
$\lamin^{x_k,y_k}\geq C y_k^{1/2}$. 
\end{step} 
 
Here and in the sequel, $C$ will always denote a
positive constant which does not depend on $k$.
 
For the proof of the first step, we let 
$\Om= \{ (x,y) \in \cM_1 \,|\, 1/2 \leq y \leq 3/2 \}$. Since $\Om$ is compact
and since  $\lamin$ is continuous and positive,
there exists $C >0$ such that for all 
\begin{eqnarray} \label{r1}
\lamin \geq C \hbox{ on } \Om.
\end{eqnarray} 
Now,
assume that 
$$\lim_k \frac{\lamin^{x_k,y_k} }{y_k^{\frac{1}{2}}} = 0.$$
 We can find a sequence $(N_k)_k$ which tends to $+\infty$
 such that $( 3^{N_k}x_k, 3^{N_k}  y_k) \in \Om$.
Note that the locally isometric covering 
$T_{p x_k, py_k}\to T_{x_k,y_k}$, $p\in \mN$, preserves the spin structures
if and only if $p$ is odd.
Let $\tilde{\psi_k}$ be the pullback of $\psi_k$ with respect to 
covering $T_{3^{N_k} x_k, 3^{N_k} y_k}\to T_{x_k,y_k}$. 
We now  have 
$$\int_{T_{3^{N_k} x_k, 3^{N_k} y_k}} {| D \psi_k |}^{\frac{4}{3}} \,dx = 
3^{N_k} \int_{T_{x_k,y_k}} {| D \psi_k |}^{\frac{4}{3}} dx$$
and  
$$\int_{T_{3^{N_k} x_k, 3^{N_k} y_k}} \big\<\psi_k ,   D {\psi_k } \big\>\, dx = 
3^{N_k}  \int_{T_{x_k,y_k}} \big\<\psi_k ,   D \psi_k \big\>\, dx.$$
We then get by (\ref{r1}) that 
\begin{eqnarray*} %\label{lamink}  
C \leq \lamin^{3^{N_k}x^k,3^{N_k}  y_k} \leq J_{3^{N_k} (x_k,y_k)}(\psi_k) =
3^{\frac{N_k}{2}}  \lamin^{x_k,y_k} \leq  C y_k^{-1/2} \lamin^{x_k,y_k}.
\end{eqnarray*}

\begin{step} \label{st3} 
There exists $C>0$ such that for all $k$,  we have $\lamin^{x_k,y_k} \geq C$. 
\end{step}

Let $\eta: \mR \to [0,1]$ be a cut-off function defined on $\mR$ 
which is equal to $0$ on $\mR
\setminus [-1,2]$ and which is equal to $1$ on $[0,1]$. We may assume that
$\eta$ is smooth. Let $v_k =(x_k,y_k)$. Since $(e_1,v_k)$ is a basis of
$\mR^2$, 
we can define $\eta_k: \mR^2 \to [0,1]$ by
$$\eta_k(tv_k+ se_1) = \eta(s)$$
Since $v_k$ is asymptotically orthogonal to $e_1$, we can find $C>0$
independent of $k$ such
that 
\begin{eqnarray} \label{gradeta}
\left|\na \eta_k\right| \leq C 
\end{eqnarray}

\noindent Moreover, by corollary \ref{cor2}, we have 
\begin{eqnarray} \label{funct} 
\frac{{\left( \int_{Z_{x_k,y_k}} {|D \eta_k \psi_k|}^{\frac{4}{3}} dx
    \right)}^{\frac{3}{2}} }{ | \int_{Z_{x_k,y_k}} \big\< \eta_k \psi_k, D
  \eta_k \psi_k\big\>
  dx|} \geq \lamin(\mS^2).
\end{eqnarray}

\noindent Now, we write that 

\begin{eqnarray*} 
{\left( \int_{Z_{x_k,y_k}} {|D \eta_k \psi_k|}^{\frac{4}{3}} dx
    \right)}^{\frac{3}{4}} & = &
{\left( \int_{Z_{x_k,y_k}} {|\nabla \eta_k \cdot \psi_k+ \eta_k D \psi_k|}^{\frac{4}{3}} dx
    \right)}^{\frac{3}{4}}\\
& \leq  & {\left( \int_{Z_{x_k,y_k}} {|\nabla \eta_k \cdot
      \psi_k|}^{\frac{4}{3}} dx 
    \right)}^{\frac{3}{4}} + {\left( \int_{Z_{x_k,y_k}}  {| \eta_kD \psi_k|}^{\frac{4}{3}} dx
    \right)}^{\frac{3}{4}}. 
\end{eqnarray*}
By (\ref{gradeta}) and H\"older inequality, we have 
\begin{eqnarray*}
{\left( \int_{Z_{x_k,y_k}} {|\nabla \eta_k \cdot  \psi_k|}^{\frac{4}{3}} dx
    \right)}^{\frac{3}{4}} & \leq  & C {\left( \int_{Z_{x_k,y_k} \cap
          \Supp(\nabla \eta_k)} {| \psi_k|}^{\frac{4}{3}} dx
    \right)}^{\frac{3}{4}} \\
& \leq & C  {\left( \int_{Z_{x_k,y_k}\cap
          \Supp(\nabla \eta_k)} 
{| \psi_k|}^4 dx
    \right)}^{\frac{1}{4}} {\Vol\Big(Z_{x_k,y_k}\cap  \Supp(\nabla \eta_k)
    \Big)}^{\frac{1}{2}}.
\end{eqnarray*}
%As one can check, $Z_{x_k,y_k} \cap
%          Supp(\eta_k)$ is contained in a parallelogramm\mnote{It is not a parallelogramm!! This phrase is more confusing than helping} which consists of
%          three copies of $T_{x_k,y_k}$ 
We then have
$$\Vol(Z_{x_k,y_k}\cap  \Supp(\nabla \eta_k) )
     \leq 3 y_k.$$
By (\ref{norm}) and step \ref{st2}, this gives that  
$${\left( \int_{Z_{x_k,y_k}} {|\nabla \eta_k \cdot  \psi_k|}^{\frac{4}{3}} dx
    \right)}^{\frac{3}{4}} \leq C y_k^{\frac{1}{2}} \leq C
  \lamin^{x_k,y_k}.$$ 
With the same argument and using relations  (\ref{eqk}) and (\ref{norm}),
it follows that 
$${\left( \int_{Z_{x_k,y_k}}  {| \eta_kD \psi_k|}^{\frac{4}{3}} dx
    \right)}^{\frac{3}{4}} \leq 3^{\frac{3}{4}} 
\lamin^{x_k,y_k} {\left( \int_{T_{x_k,y_k}}
{| \psi_k|}^4 dx
    \right)}^{\frac{3}{4}} \leq C \lamin^{x_k,y_k}.$$
Finally, we get that  
\begin{eqnarray}  \label{numer1}
{\left( \int_{Z_{x_k,y_k}} {|D \eta_k \psi_k|}^{\frac{4}{3}} dx
     \right)}^{\frac{3}{2}} \leq C  {\left(\lamin^{x_k,y_k} \right)}^2.
\end{eqnarray}
We now write that 
$$ \int_{Z_{x_k,y_k}} \< \eta_k \psi_k, D
  \eta_k \psi_k\>\,dx = 
\int_{Z_{x_k,y_k}} \< \eta_k \psi_k,\nabla \eta_k \cdot \psi_k + \eta_k D
   \psi_k\>\, dx.$$
Moreover, the left hand side of this equality is real since $D$ is an 
autoadjoint operator. Since 
$$\int_{Z_{x_k,y_k}} \< \eta_k \psi_k,\nabla \eta_k \cdot \psi_k\>\, dx
  \in i \mR.$$ 
Together with equation (\ref{eqk}), this implies that 
$$ \int_{Z_{x_k,y_k}} \big\< \eta_k \psi_k, D
  \eta_k \psi_k\big\>
  dx  =  \int_{Z_{x_k,y_k}} \eta_k^2 \lamin^{x_k,y_k} |\psi_k|^4 dx.$$
Using (\ref{norm}), we obtain that 
\begin{eqnarray} \label{denom1} 
 \int_{Z_{x_k,y_k}} \big\< \eta_k \psi_k, D
  \eta_k \psi_k\big\>\,
  dx \geq  \lamin^{x_k,y_k} \int_{T_{x_k,y_k}}   |\psi_k|^4 dx =
  \lamin^{x_k,y_k}.
\end{eqnarray}
Finally, plugging  (\ref{numer1}) and (\ref{denom1}) in (\ref{funct}), we
obtain that $\lamin(\mS^2) \leq C \lamin^{x_k,y_k}$.
This proves the step.

\begin{step} \label{st4} 
The function $\lamin$  can be extended continuously to  
$\cM_1 \cup \{(0,0)\}$ by  setting $ \lamin^{0,0} =  \lamin(\mS^2)$. 
\end{step} 
  
\noindent In other words, we show that $\lim_k  \lamin^{x_k,y_k} =  \lamin(\mS^2)$.
The method is quite similar than the one of previous step.
 Let $\zeta_k: \mR \to [0,1]$ be a smooth cut-off function defined on $\mR$ 
which is equal to $0$ on $\mR
\setminus [-y_k,1+y_k]$, which is equal to $1$ on $[0,1]$ and which satisfies
$|\nabla \zeta_k| \leq \frac{2}{y_k}$. As in the last step, we can define
 $\gamma_k: \mR^2 \to [0,1]$ by
$$\gamma_k(tv_k+ se_1) = \zeta_k(s).$$
Since $v_k$ is asymptotically orthogonal to $e_1$, we can find $C>0$
independent of $k$ such
that 
\begin{eqnarray} \label{gradgamma}
|\gamma_k| \leq \frac{C}{y_k} 
\end{eqnarray}
As in step \ref{st3}, we have 
\begin{eqnarray} \label{funct1} 
\frac{{\left( \int_{Z_{x_k,y_k}} {|D \gamma_k \psi_k|}^{\frac{4}{3}}\, dx
    \right)}^{\frac{3}{2}} }{|  \int_{Z_{x_k,y_k}} \big\< \gamma_k \psi_k, 
    D  \gamma_k \psi_k\big\>\,dx|} \geq \lamin(\mS^2).
\end{eqnarray}

\noindent We first prove that we can assume that 
\begin{eqnarray} \label{supp}
  \int_{Z_{x_k,y_k} \cap \Supp(\nabla \gamma_k)}   |\psi_k|^4 dx  \leq
  {C}{y_k}.
\end{eqnarray} 
We let $n_k= [{(2y_k)}^{-1}]$ be the integer part of ${2 y_k}^{-1}$. For all 
$l \in [0,n_k-1]$, we define 
$$A_{k,l}= \Big\{ t e_1+ sv_k | s \in [0,1[ \hbox{ and } t \in
\left[ \frac{l-\frac{1}{2}}{n_k}, \frac{l+\frac{1}{2}}{n_k} \right] \Big\}.$$
The family of sets $(A_{k,l})_{l \in [0,n_k-1]}$ is a partition of 
$T'_{x_k,y_k}$ which is the image of $T_{x_k,y_k}$ by the translation of
vector $-\frac{1}{2n_k}e_1 $. By periodicity, $(A_{k,l})_{l \in [0,n_k-1]}$
can be seen as
a partition of $T_{x_k,y_k}$.
Consequently, we can write that 
$$1= \int_{T_{x_k,y_k}} {|\psi_k|}^4 dx= \sum_{l=0}^{n_k-1} \int_{A_{k,l}}
{|\psi_k|}^4 dx.$$
Hence, there exists $l_0  \in [0,n_k-1]$ such that 
$$ \int_{A_{k,l_0}}
{|\psi_k|}^4 dx = \min_{l \in [0,n_k-1]} \sum_{l=0}^{n_k-1} \int_{A_{k,l}}
{|\psi_k|}^4 dx \leq \frac{1}{n_k}.$$
Obviously, without loss of generality, we can replace $\psi_k$ by $\psi_k
\circ t_0$ where $t_0$ is the translation of vector $-l_0 e_1$. In this
way, we can assume that $l_0=0$. By periodicity,
$\Supp(\nabla \gamma_k) \subset  
A_{k,0}$. Hence,
$$\int_{Z_{x_k,y_k} \cap \Supp(\nabla \gamma_k)}   |\psi_k|^4 dx  \leq
  \frac{1}{n_k}.$$
Since $n_k \sim \frac{2}{y_k}$, equation (\ref{supp}) follows.

\noindent Now, we proceed as in step \ref{st3}. We write that

\begin{eqnarray*}
{\left( \int_{Z_{x_k,y_k}} {|D \gamma_k \psi_k|}^{\frac{4}{3}} dx
    \right)}^{\frac{3}{4}} &  =  &
{\left( \int_{Z_{x_k,y_k}} {|\nabla \gamma_k \cdot \psi_k+ \gamma_k D \psi_k|}^{\frac{4}{3}} dx
    \right)}^{\frac{3}{4}}\\
& \leq & {\left( \int_{Z_{x_k,y_k}} {|\nabla \gamma_k \cdot
      \psi_k|}^{\frac{4}{3}} dx 
    \right)}^{\frac{3}{4}} + {\left( \int_{Z_{x_k,y_k}}  
 {| \gamma_k D \psi_k|}^{\frac{4}{3}} dx
    \right)}^{\frac{3}{4}}. 
\end{eqnarray*}
It follows from (\ref{gradgamma}) and the H\"older inequality that 
\begin{eqnarray*}
 {\left( \int_{Z_{x_k,y_k}} {|\nabla \gamma_k \cdot  \psi_k|}^{\frac{4}{3}} dx
    \right)}^{\frac{3}{4}} & \leq & \frac{C}{y_k} 
 {\left( \int_{Z_{x_k,y_k} \cap
          \Supp(\nabla \gamma_k)} {| \psi_k|}^{\frac{4}{3}} dx
    \right)}^{\frac{3}{4}}\\
& \leq &\frac{C}{y_k}   {\left( \int_{Z_{x_k,y_k}\cap
          \Supp(\nabla \gamma_k)} 
{| \psi_k|}^4 dx
    \right)}^{\frac{1}{4}} {\Big( \Vol(Z_{x_k,y_k}\cap  \Supp(\nabla \gamma_k) )
    \Big)}^{\frac{1}{2}}.
\end{eqnarray*}
Clearly, we have 
$$\Vol(Z_{x_k,y_k}\cap  \Supp(\nabla \gamma_k) )
     \leq C y_k^2.$$
By (\ref{supp}), we obtain  
$${\left( \int_{Z_{x_k,y_k}} {|\nabla \gamma_k \cdot  \psi_k|}^{\frac{4}{3}} dx
    \right)}^{\frac{3}{4}} \leq C y_k^{-1+ \frac{1}{4}+1} \leq C 
 y_k^{\frac{1}{4}}= o(1). $$ 
For the other term, we write, using (\ref{eqk})
$${\left( \int_{Z_{x_k,y_k}}  {| \gamma_k D \psi_k|}^{\frac{4}{3}} dx
    \right)}^{\frac{3}{4}} =  
\lamin^{x_k,y_k} {\left(  \int_{T_{x_k,y_k}}
{| \psi_k|}^4 dx +  \int_{Z_{x_k,y_k} \cap \{ 0< \gamma_k <1 \} } 
{| \psi_k|}^4 dx
    \right)}^{\frac{3}{4}}.$$
Clearly, we can construct $\gamma_k$  such that $ \{ 0< \gamma_k <1 \}
\subset \Supp( \nabla \gamma_k) $. It then follows from (\ref{supp}) that 
$${\left( \int_{Z_{x_k,y_k}}  {| \gamma_k D \psi_k|}^{\frac{4}{3}} dx
    \right)}^{\frac{3}{4}} \leq \lamin^{x_k,y_k} + o(1).$$
Finally, we obtain 
\begin{eqnarray}  \label{numer2}
{\left( \int_{Z_{x_k,y_k}} {|D \gamma_k \psi_k|}^{\frac{4}{3}} dx
     \right)}^{\frac{3}{2}} \leq   {\left(\lamin^{x_k,y_k} \right)}^2+ o(1).
\end{eqnarray}
Now, as in step \ref{st3}, we write that 
$$ \int_{Z_{x_k,y_k}} \big\< \gamma_k \psi_k, D
  \gamma_k \psi_k\big\>\,
  dx 
 =  \int_{Z_{x_k,y_k}} \gamma_k^2 \lamin^{x_k,y_k} |\psi_k|^4 dx.$$
Using (\ref{norm}), we obtain that 
\begin{eqnarray} \label{denom2} 
 \int_{Z_{x_k,y_k}} \big\< \gamma_k \psi_k, D
  \gamma_k \psi_k \>\,
  dx \geq  \lamin^{x_k,y_k} \int_{T_{x_k,y_k}}   |\psi_k|^4 dx =
  \lamin^{x_k,y_k}.
\end{eqnarray}
Plugging  \eref{numer2} and \eref{denom2} in \eref{funct1}, we
obtain that 
  $$\lamin(\mS^2) \leq  {\left(\lamin^{x_k,y_k}\right)^2 + o(1)\over 
    \lamin^{x_k,y_k}}$$ 
which implies that either
$\lamin^{x_k,y_k}\to 0$ or $\lamin^{x_k,y_k}\to \lamin(\mS^2)$. 
Hence, step \ref{st3} yields the statement of the theorem.
\qed 
 
%\section{Proof of corollary \ref{cor}} 
%At first, assume that  $n \geq 2$. Let $g$ be a Riemannian  
%metric on $\mS^1 \times \mS^n$. We let $\mu([g])$ be the Yamabe invariant   
%of $(\mS^1\times \mS^n,g)$. According to \cite{}, we have 
%$$\sup \mu([g]) = \mu(\mS^n)$$ 
%where the supremum is taken  over the set of conformal classes of metrics on  
%$\mS^1 \times \mS^n$. 
%Moreover, Hijazi inequality asserts that  
%$$ \lamin([g])  \geq \frac{n}{4(n-1)}  \mu([g]) $$ 
%It is well known that  
%$$\lamin(\mS^n) = \frac{n}{4(n-1)}  \mu(\mS^n)$$ 
%Corollary \ref{cor} then follows. 
 
%\noindent We assume now that $n=1$ and we note $T=\mS^1 \times \mS^1 = [0,2
%\pi[ \times [0,2 \pi[$. 
%%Let $(e_1,e_2)$ be the canonical basis of  
%%$\mR^2$. We let $v=(x_0,y_0) \in \mR^2$ then  
 
%%\[ \left| \begin{array}{ccc}  (\T,g_{eucl})  & \to & (T, g_{x_0,y_0})\\ 
%%x e_1 + x v & \mapsto & 2 \pi (x e_1 + x e_2) \end{array}  
%%\right. \] 
%%where $g_{x_0,y_0}$ is the  metric on $T$ defined in $(e_1,e_2)$ by  
%%$$ \frac{1}{2 \pi} \left(\begin{matrix}  
%%1 & (e_2,v)\\ 
%%(e_2,v)  & x_0^2 + y_0^2 
%%\end{matrix} \right) $$ 
%%Here, $(.,.)$ denotes the canonical scalar product on $\mR^2$. 
%\noindent By lemma \ref{isom}, to each $(x_0,y_0) \in \cM_1$ corresponds a metric
%on the torus $T=\mS^1 \times \mS^1$. Moreover, by theorem \ref{main}   
%$$ \lamin(\mS^2) = sup_{ (x_0,y_0) \in \cM_1} \lam 
%\leq \sup_{g \hbox{ metric on T}} \lamin([g])$$ 
%This proves corollary (\ref{cor}). 

%\bibliographystyle{amsbernd}
%\bibliography{literatur}

\begin{thebibliography}{AHM60}            
          
%\bibitem[AAF99]{agricola.ammann.friedrich:99}
%I.~Agricola, B.~Ammann, T.~Friedrich,
%\emph{A comparison of the eigenvalues of the Dirac and Laplace operator 
%on the two-dimensional torus},
%Manus. Math. {\bf 100} (1999), 231--258.

\bibitem[Am98]{ammann:diss}
B.~Ammann,
\newblock {\em {S}pin-{S}trukturen und das {S}pektrum des {D}irac-{O}perators},
\newblock PhD thesis, University of {F}reiburg, {G}ermany, 1998.
\newblock Shaker-Verlag Aachen 1998, ISBN 3-8265-4282-7.

\bibitem[Am03]{ammann:03}    
B.~Ammann, \emph{A spin-conformal lower bound of the first positive {D}irac    
  eigenvalue}, {D}iff. {G}eom. {A}ppl. \textbf{18} (2003), 21--32.    

\bibitem[Am03a]{ammann:habil}    
B.~Ammann, \emph{A variational problem in conformal spin geometry},  
Habilitationsschrift, Universit\"at Hamburg, May 2003,   
dowloadable on http://www.berndammann.de/publications.  

\bibitem[Am03b]{ammann:p03}
B.~Ammann, \emph{The smallest {D}irac eigenvalue in a spin-conformal class and
  cmc-immersions}, {P}reprint, 2003.

\bibitem[AB02]{ammann.baer:02}B.~Ammann, C.~B\"ar, 
\emph{Dirac eigenvalue estimates on surfaces}, 
Math. Z. {\bf 240} (2002), 423--449. 


%\bibitem[AH03]{ammann.humbert:03}          
%B.~Ammann, E. Humbert,    
%\newblock{\em Positive mass theorem for the Yamabe problem on spin manifolds}, \newblock{Preprint}, to appear
%in Geom. Funct. Anal.            

%\bibitem[AHM03]{ammann.humbert.morel:p03b}          
%B.~Ammann, E. Humbert, B. Morel,   
%\newblock{\em Mass endomorphism and spinorial Yamabe type problems   
%on conformally flat manifolds.}   
%\newblock{Preprint}            
                  
\bibitem[AHM03]{ammann.humbert.morel:03}          
B.~Ammann, E. Humbert, B. Morel,   
\newblock{\em A spinorial analogue of Aubin's inequality,}   
\newblock{Preprint},  ArXiv math.DG/0308107. 

\bibitem[Aub98]{aubin:98}
T.~Aubin,
\newblock \emph{Some Nonlinear problems in Riemannian Geometry},     
\newblock {Springer-Verlag, 1998. }

    
\bibitem[B\"ar92]{baer:92b}            
C.~B\"ar,            
\newblock \emph{Lower eigenvalue estimates for Dirac operators},            
\newblock {Math. Ann.}, {\bf 293} (1992), 39--46.            
     
\bibitem[B\"ar96]{baer:96}C.~B\"ar,            
\emph{Metrics with harmonic spinors}, 
Geom.\ Funct.\ Anal.\ {\bf 6} (1996), 899--942.

\bibitem[BD02]{baer.dahl:02}
C.~B{\"a}r and M.~Dahl, \emph{{Surgery and the Spectrum of the {D}irac
  Operator}}, J. reine angew. Math. \textbf{552} (2002), 53--76.

\bibitem[BS92]{baer.schmutz:92}
C.~B{\"a}r and P.~Schmutz, \emph{{Harmonic spinors on Riemann surfaces}}, Ann.
  Global Anal. Geom. \textbf{10} (1992), 263--273.

\bibitem[BGV91]{berline.getzler.vergne:91}
N. Berline and E. Getzler and M. Vergne,
\emph{Heat kernels and {D}irac operators},
Springer Verlag 1991.


%\bibitem[BG92]{bourguignon.gauduchon:92}  
%J.-P.~Bourguignon and P. Gauduchon,  
%\newblock {\em Spineurs, op{\'e}rateurs de {D}irac et variations de m{\'e}triques},  
%\newblock {Comm. Math. Phys.} \textbf{144} (1992), 581--599.  
%

\bibitem[BN04]{bray.neves:04}
H.~L. Bray and A.~Neves, \emph{Classification of prime 3-manifolds with
  {Y}amabe invariant greater than {$\mathbb{R}P\sp 3$}} Ann. of Math.
  \textbf{159} (2004), 407--424.
 
%\bibitem[BN04]{bray.neves:04a}
%H.~L. Bray and A.~Neves, \emph{Corrigendum: ``{C}lassification of prime
%  3-manifolds with {Y}amabe invariant greater than {$\mathbb{R}P\sp 3$}''}, 
%Ann. of Math. (2) \textbf{159} (2004), 887.

\bibitem[Fr84]{friedrich:84}
T.~Friedrich,
{\it Zur Abh\"angigkeit des Dirac-Operators von der Spin-Struk\-tur},
Colloq.~Math. {\bf 48} {\it 57--62} (1984)

\bibitem[Fr00]{friedrich:00}
T.~Friedrich,
\emph{{D}irac {O}perators in {R}iemannian {G}eometry},
{{G}raduate {S}tudies in {M}athematics 25},
{AMS, Providence, Rhode Island},
2000

\bibitem[GL80]{gromov.lawson:80}
M.~Gromov and H.~B. Lawson, \emph{The classification of simply connected
  manifolds of positive scalar curvature}, Ann. of Math. \textbf{111}
  (1980), 423--434.
%
%\bibitem[Heb97]{hebey:97}            
%E. Hebey,            
%\newblock  {\em Introduction \`a l'analyse non-lin\'eaire sur les vari\'et\'es.}    
%\newblock {Diderot \'Editeur, Arts et sciences.}, 1997.         
    
\bibitem[Hij86]{hijazi:86}            
O.~Hijazi,            
  \newblock  \emph{A conformal lower bound for the smallest eigenvalue of the {D}irac operator and {K}illing Spinors}, 
\newblock {Comm. Math. Phys.} \textbf{104} (1986), 151--162.  
            
\bibitem[Hij91]{hijazi:91}              
O.~Hijazi,              
\newblock  \emph{Premi\`ere valeur propre de l'op\'erateur de {D}irac et nombre de {Y}amabe,}              
\newblock {C. R. Acad. Sci. Paris}, \newblock{\em S\'erie I } \textbf{313}, (1991), 865--868.             

\bibitem[Hij01]{hijazi:01}
O.~Hijazi,              
\emph{Spectral properties of the {D}irac operator and geometrical structures},
Ocampo, Hernan (ed.) et al., Geometric methods for quantum field theory. 
Proceedings of the summer school, Villa de Leyva, Colombia, July 12-30, 1999. 
Singapore: World Scientific. 116-169,
2001.
  
\bibitem[Hit74]{hitchin:74}
N.~Hitchin, \emph{Harmonic spinors}, Adv. Math. \textbf{14} (1974), 1--55.

%\bibitem[LP87]{lee.parker:87}            
%J.~M. Lee and T.~H. Parker.            
%\newblock \emph{The Yamabe problem.}            
%\newblock {Bull. Am. Math. Soc., New Ser.} \textbf{17} (1987), 37--91.            

\bibitem[LM89]{lawson.michelsohn:89}
H.-B.~Lawson and M.-L. Michelsohn.
\newblock {\em Spin Geometry}.
\newblock Princeton University Press, Princeton, 1989.

\bibitem[Lott86]{lott:86}            
J.~Lott,            
\newblock \emph{Eigenvalue bounds for the Dirac operator,}            
\newblock {Pacific J. of Math.} \textbf{125} (1986), 117--126.              

\bibitem[Mai97]{maier:97}
S.~Maier, \emph{Generic metrics and connections on spin- and
  spin-$\,^c$-manifolds}, Comm. Math. Phys. \textbf{188} (1997), 407--437.

\bibitem[Mil63]{milnor:63}
J.~Milnor, \emph{Spin structures on manifolds}, Enseignement Math. (2)
  \textbf{9} (1963), 198--203.

\bibitem[Ro88]{roe:88}
J. Roe,
\emph{Elliptic operators, topology and asymptotic methods},
{Pitman Research Notes in Mathematics Series},
{\bf 179}, Longman 1988.
            
%\bibitem[Sch84]{schoen:84}            
%R. Schoen.            
%\newblock {Conformal deformation of a Riemannian metric to constant scalar            
%  curvature}.            
%\newblock {\em J. Diff. Geom.}, \textbf{20} (1984), 479--495.            
            
%\bibitem[Sch91]{schoen:91}            
%R.~Schoen,            
%\newblock \emph{On the number of constant scalar curvature metrics in a
%  conformal class}, 
%\newblock {Differential Geometry: A symposium in honor of Manfredo Do
%  Carmo}, Pitman Monog. Surveys Pure Appl. Math., \textbf{54},  Longman Sci. Tech., Harlow (1984), 
%311--320.            

\bibitem[Sch91]{schoen:91}            
R.~Schoen, \emph{On the number of constant scalar curvature metrics in a
conformal class},
\newblock {Differential Geometry: A symposium in honor of Manfredo Do
Carmo (B. Lawson et K. Teneblat (Ed.)},
{Pitman Monogr. Surveys Pure Appl. Math.}, \textbf{52},
Longman Sci. Tech., Harlow (1991), 311--320.


\bibitem[St92]{stolz:92}
S.~Stolz, \emph{Simply connected manifolds of positive scalar curvature},
Ann. of Math. (2),
 {\bf 136} (1992), {\it 511--540}.

\bibitem[Ta81]{taylor:81}
M.~E.~Taylor,
\emph{Pseudodifferential operators},
Princeton University Press, 1981.
            
\end{thebibliography}

\vspace{1cm}              
Authors' address:              
\nopagebreak  
\vspace{5mm}\\  
\parskip0ex    
\vtop{  
\hsize=8cm\noindent  
\obeylines              
Bernd Ammann and Emmanuel Humbert              
Institut \'Elie Cartan BP 239             
Universit\'e Henri Poincar\'e, Nancy              
54506 Vandoeuvre-l\`es -Nancy Cedex              
France                          
\vspace{0.5cm}              
              
E-Mail:              
{\tt ammann at iecn.u-nancy.fr and humbert at iecn.u-nancy.fr}              
}  
                   
%%%%%%%%%%%%%%%%%%%%%%%%%%%%%%%%%%%%%%%%%%%%%%%%%%%%%%%%%%%%%%%%%    
\end{document}